\documentclass[journal]{IEEEtran}

\usepackage{cite}

\ifCLASSINFOpdf
   \usepackage[pdftex]{graphicx}
\else
\fi

\usepackage{amsmath}

\begin{document}

\title{Distribution Network Marginal Costs --- Part I: A Novel AC OPF Including Transformer Degradation}

\author{Panagiotis~Andrianesis,
        and~Michael~Caramanis,~\IEEEmembership{Senior~Member,~IEEE}
\thanks{P. Andrianesis and M. Caramanis are with the Division of Systems Engineering, Boston University, 
Boston, MA, 02446 USA, e-mails: panosa@bu.edu, mcaraman@bu.edu. Research partially supported by the Sloan Foundation under grant G-2017-9723 and NSF AitF grant 1733827.}
}

\maketitle

\begin{abstract}
This two-part paper considers the day-ahead operational planning problem of a radial distribution network hosting Distributed Energy Resources (DERs) including Solar Photovoltaic (PV) and storage-like loads such as Electric Vehicles (EVs). 
We estimate dynamic Distribution nodal Location Marginal Costs (DLMCs) of real and reactive power and employ them to co-optimize distribution network and DER schedules. 
In Part I, we develop a novel AC Optimal Power Flow (OPF) model encompassing transformer degradation as a short-run network variable cost, and we decompose real/reactive power DLMCs into additive marginal cost components related to (\emph{i}) the costs of real/reactive power transactions at the T\&D interface/substation, (\emph{ii}) real/reactive power marginal losses, (\emph{iii}) voltage and (\emph{iv}) ampacity congestion, and (\emph{v}) a new transformer degradation marginal cost component.
Our detailed transformer degradation model captures the impact of incremental transformer loading during a specific time period, not only on its Loss of Life (LoL) during that period, but also during subsequent time periods. 
To deal with this phenomenon, we develop methods that internalize the marginal LoL occurring beyond the daily horizon into the DLMCs evaluated within this horizon.
In Part II, we use real distribution feeders to exemplify the use of DLMCs as financial incentives that convey sufficient information to optimize Distribution Network, and DER (PV and EV) operation.
\end{abstract}


\IEEEpeerreviewmaketitle

\section{Introduction}

\IEEEPARstart{R}{apidly} growing adoption of Distributed Energy Resources (DERs), including clean, albeit volatile, renewable generation, combined heat and power micro generation, storage and flexible loads with storage-like properties and Volt/VAR control capabilities, e.g., Electric Vehicles (EVs) and Solar Photovoltaic (PV) inverters, present a major challenge together with still unexploited opportunities.
As DERs become a major user of distribution grid infrastructure, the grid is becoming increasingly active, distributed, dynamic, and challenging to plan and operate \cite{Tierney2017}.
DERs are bound to have a profound impact on the adequacy of T\&D assets, efficient grid operation, reliability, and security of supply, and as such their scheduling will be crucial. 
This paper shows that optimal DER scheduling depends on their dynamic Distribution Locational Marginal Cost (DLMC) \cite{CaramanisEtAl2016, TaborsEtAl2016, NtakouEtAlIEEETD2014} whose accurate estimation is expected to bring about fundamental changes in distribution planning, operation, and eventually power markets.

\subsection{Background and Motivation}
There is a rich literature on DER integration into Distribution Networks. 
Numerous studies address the impact of various DERs (primarily EVs and PVs) on the grid and its assets ---  e.g., overvoltages due to PV --- and examine DER capabilities --- e.g., the provision of reactive power \cite{LopesEtAl_2007, OlivierEtAl_2016, GongEtAl_2012, HilsheyEtAl_2013, QianEtAl_2015, AssolamiMorsi2015, ElNozahySalama_2014, PezeshkiEtAl_2014, GrayMorsi_2017, WangEtAl_2018}. 
Several works investigate also the impact of EVs on distribution transformers \cite{GongEtAl_2012, HilsheyEtAl_2013, QianEtAl_2015, AssolamiMorsi2015, ElNozahySalama_2014}.
They focus on the acceleration of transformer degradation from increasing EV penetration, noting that transformer aging is dependent upon the thermal effects of persistent transformer loading.

The life of a transformer is strongly related to its winding hottest spot temperature (HST) driving insulation aging.
Distribution transformer life may exceed 20 years (180,000 hours), assuming operation at a reference HST of $110^o$C (insulation designed for an average winding temperature of $65^o$C).
IEEE Standard C.57.91-2011 \cite{IEEE2011} and IEC Standard 60076-7:2018 \cite{IEC2018} provide guidelines for transformer loading and propose an exponential representation of the aging acceleration factor in terms of the HST when it exceeds $110^o$C.
In particular, for HSTs exceeding (being lower than) $110^o$C, the aging acceleration factor is greater (less) than 1.

Early EV pilot studies (e.g., \cite{EVproject2013}) identified that clustering EV chargers under the same transformer may cause damage and outages from persistent overloading.
Recently, a study for the Sacramento Municipal Utility District \cite{SMUD} estimated that, under a high DER penetration scenario, about 26\% of the substations and service transformers would experience overvoltages (by at least 5\% of nominal) due to PVs, whereas up to 17\% of the approximately 12,000 service transformers would experience overloads (exceeding 140\% the nameplate rating) due to EVs and would need to be replaced at an average estimated cost of \$7,400 per transformer.
The 2017 FERC Financial Form I filing by Commonwealth Edison, a typical urban distribution Utility, reports the original cost of line transformers at 8\% of the total cost of the Electric Plant. 

Unsurprisingly, DER scheduling and coordination with the grid and its assets in the context of distribution network operational planning is attracting increasing attention. 
Emerging literature focuses on the extension of wholesale Locational Marginal Price (LMP) markets to the distribution network through reliance on Distribution LMPs (DLMPs).
A variety of suggested approaches and models have attempted to do just that by considering uniform price-quantity bidding DERs.
This two-part paper argues that most published work is lacking in considering the daily cycle costs of distribution networks and the complex preferences of DERs which are inadequately represented by uniform price quantity bids and offers. Indeed, past work has paid insufficient attention to Network-DER co-optimization and the role of economically efficient DLMCs that go beyond wholesale market LMPs by internalizing \emph{intertemporally} coupled DER preferences (e.g., energy rather than capacity demanding EVs), and salient distribution network costs (e.g., loss of transformer life resulting from persistent loading).

The considerable time that elapsed between the proposal of wholesale power markets in the eighties \cite{CaramanisEtAl_1982, SchweppeEtAl_1988} and their eventual, and moreover extensive and successful, implementation in the late nineties, behooves us to redouble our effort in the appropriate modeling and estimation of DLMCs.
The expected shock that distributed renewable generation and mobile EV battery charging will most likely deliver to distribution assets (transformers in particular) may dwarf the shock delivered in the late seventies by the massive adoption of air-conditioning. 
What costs should DLMCs capture to promote efficient scheduling of EV charging and dual use of PV smart inverters?
Trading of reactive power was considered but abandoned since it was deemed of limited value in the wholesale transmission markets designed in 1997.
However, nodal voltage congestion and sustained transformer overloading in distribution networks render the explicit modeling of reactive power a \emph{must} in the context of a 24-hour ahead AC OPF capable of modeling Volt/VAR control, ampacity constrained feeders, intertemporally coupled transformer life degradation and complex DER/EV charging preferences.

\subsection{Literature Review}

Existing literature investigating the EV impact on distribution transformers \cite{GongEtAl_2012, HilsheyEtAl_2013, QianEtAl_2015, AssolamiMorsi2015, ElNozahySalama_2014} includes mostly quantitative simulation studies considering various charging schedules and assessing their impact on transformer Loss of Life (LoL).
Indicative results show that the daily LoL of a residential service transformer may almost double if EVs charge upon arrival \cite{QianEtAl_2015}, and that simple ``open-loop'' EV charging, e.g., delay until after midnight that may be promoted by Time of Use rates, can actually increase, rather than decrease, transformer aging under likely circumstances \cite{HilsheyEtAl_2013}. 

Rooftop PVs, on the other hand, although they may create overvoltage issues, they can utilize their inverters to mitigate overvoltages and most importantly to schedule their real and reactive power output with a beneficial impact on transformer LoL.
Indeed, findings \cite{PezeshkiEtAl_2014} show that high PV penetration can significantly extend the life of distribution transformers in a suburban area.
Furthermore, while previous research determined insignificant synergy between the substation transformer LoL and charging EVs in the presence of rooftop PV (due to non-coincidence between peak hours of PV generation and EV charging), \cite{GrayMorsi_2017} shows that the transformer thermal time constant allows PV generation to reduce the transformer temperature when EVs are charging.

Although the above works consider the DER impact on transformer aging, they do not internalize transformer degradation to the operational planning problem.
A first outline is provided in \cite{CaramanisEtAl2016, TaborsEtAl2016, NtakouEtAlIEEETD2014} --- among the first DLMP works --- in which the exponential transformer LoL model is employed in an AC OPF formulation with nonlinear LoL objective function components but omission of the intertemporal coupling in transformer LoL associated with heating and cooling time constants.
In \cite{CaramanisEtAl2016}, DLMPs are extended to include reserves --- apart from real/reactive power --- and an iterative distributed architecture is sketched capturing DER intertemporal preferences and physical system dynamics.
These works employ the \emph{branch} \emph{flow} (a.k.a. DistFlow) equations, introduced for radial networks by \cite{BaranWu1989}, and in fact their relaxation to convex Second Order Cone constraints proposed by \cite{FarivarLow2013}.

Parallel works \cite{LiEtAl_2014, HuangEtAl_2015} determined DLMPs in a social welfare optimization problem, using DC OPF, and considering EV aggregators as price takers.
Quadratic programming is used in \cite{HuangEtAl_2015}, as opposed to linear programming used in \cite{LiEtAl_2014}, to derive DLMPs that are announced to aggregators of EVs and heat pumps optimizing their energy plans. 

Following \cite{CaramanisEtAl2016}, \cite{BaiEtAl-DLMP} presents a detailed AC OPF formulation, also employing the relaxed branch flow model, acknowledging that reactive power and voltage are critical features in the operation of distribution networks, and, hence, DC OPF approaches are inadequate in capturing salient features of distribution networks.
It proposes a day-ahead market-clearing framework for radial networks, in which various DERs, such as distributed generators, microgrids, and load aggregators bid into a day-ahead distribution-level electricity market.
DLMPs are derived using sensitivity factors of the linear version of the DistFlow model;
they include ampacity congestion but do not consider the impact on transformers.

In a parallel to \cite{BaiEtAl-DLMP}, \cite{Papavasiliou_2017} provides a comprehensive analysis of various approaches in decomposing and interpreting DLMPs in the AC OPF context, and discusses accuracy and computational aspects of the DLMP decomposition.
It considers real/reactive power and losses, as well as voltage and line congestion (modeling power instead of ampacity limits).
It does not include transformer degradation.
The analysis presents three approaches that are based on: 
(\emph{i}) duality of the Second Order Cone Programming (SOCP) formulation, 
(\emph{ii}) an implicit function formulation following \cite{CaramanisEtAl2016}, and
(\emph{iii}) an interpretation based on marginal losses.
For the latter two approaches, which involve the calculation of partial derivatives of primal variables \emph{w.r.t.} the net demand, \cite{Papavasiliou_2017} suggests the solution of several slightly perturbed AC OPF problems and shows that, apart from the increased computational requirements compared to the SOCP duality analysis, this numerical approximation may introduce errors in the calculation of the DLMP components.


\subsection{Objective and Contribution} \label{Contribution}

In this work, our main objective is to derive the time and location specific marginal costs of real and reactive power in the distribution network that are consistent with the optimal DER schedule.
As such, we aim at discovering the DLMCs and understanding the components, building blocks, and sources that constitute the marginal costs.
The main contributions of Part I of this two-part paper are as follows.

First, we present an enhanced AC OPF operational planning formulation, which includes a detailed cost representation of transformer degradation in the presence of DERs (PVs, EVs).
More specifically, we enhance the relaxed SOCP branch flow problem for radial networks, with accurate modeling of the distribution transformer LoL that takes into account its intertemporal thermal response.
In fact, we present an exact relaxation of the transformer LoL, which exploits the features of the branch flow model.
This is the first SOCP AC OPF formulation of the operational planning problem that considers a full transformer degradation model and PV/EV scheduling.
For completeness, we also provide a concise representation of the PV/EV capabilities.
Specifically, we propose general formulations that enable the representation of smart inverter capabilities and EV mobility in a multi-period operational planning problem.

Second, we provide intuitive formulas of the real and reactive power DLMC components that are related to the transformer degradation cost.
More specifically, we employ sensitivity analysis to decompose DLMCs into additive cost components, that relate to the cost of real and reactive power at the substation, the cost of real and reactive power marginal losses, voltage congestion, line (ampacity) congestion, and transformer degradation.
We illustrate how the latter component informs us on the intertemporal impact of the thermal response dynamics on the marginal cost.
We quantify this impact showing how the transformer load at a specific hour affects subsequent hours and how it reflects on the marginal degradation cost component.

Third, we discuss practical considerations for capturing the transformer degradation impact that extends beyond the day-ahead horizon and we quantify this impact on the transformer degradation DLMC component.

Lastly, we suggest an improvement compared to the analysis provided in \cite{Papavasiliou_2017}, illustrating that the partial derivatives required to calculate the marginal cost components can be derived by the solution of a linear system, thus reducing the computational effort and removing numerical approximation error concerns associated with the solution of perturbed AC OPF problems.
We also note that our analysis and derivations of the marginal cost components are different compared to \cite{BaiEtAl-DLMP}, which uses the linear version of the branch flow model.

In Part II \cite{PartII} of this two-part paper, we employ DLMCs as price signals that provide sufficient information for DERs to self-schedule in an fashion that optimizes network and DER costs.
We also explore a realistic distribution feeder and employ various inelastic demand and DER penetration scenarios to illustrate the superiority of our approach relative to reasonable conventional scheduling alternatives.
Most importantly, the realistic numerical results show the paramount importance of modeling the intertemporally coupled transformer LoL, which achieves variable network asset cost reductions that far exceed the benefits of merely enforcing voltage constraints. 

\subsection{Paper Organization}
The remainder of this paper is organized as follows.
Section \ref{Model} provides the detailed model formulation of the AC OPF operational planning optimization problem.
Section \ref{DLMCs} presents the decomposition of DLMCs with emphasis on the new transformer degradation component.
Section \ref{Conclusions} concludes and provides directions for further research.
The details of transformer HST calculation are included in Appendix \ref{AppA}.

\section{Enhanced AC OPF Model} \label{Model}
In this section, we present the formulation of an enhanced AC OPF model of the operational planning problem in radial distribution networks.
We introduce the network model and notation in Subsection \ref{NetworkModel}.
We provide a concise representation of DER (PV/EV) models in Subsection \ref{DERconstraints}, and the detailed transformer degradation formulation 
in Subsection \ref{TransformerModel}.
For ease of exposition, HST calculations are listed in Appendix \ref{AppA}.
We summarize the optimization problem in Subsection \ref{ProblemSummary}.

\subsection{Network Model} \label{NetworkModel}

We consider a radial network with $N+1$ nodes and $N$ lines.
Let $\mathcal{N}= \{ 0,1,...,N \}$ be the set of nodes, with node $0$ representing the root node, and $\mathcal{N^+} \equiv \mathcal{N} \setminus \{0\}$.
Let $\mathcal{L}$ be the set of lines, with each line denoted by the pair of nodes $(i,j)$ it connects --- henceforth $ij$ for short, where node $i$ refers to the (unique due to the radial structure) preceding node of $j \in \mathcal{N^+}$.
Transformers are represented as a subset of lines, denoted by $y \in \mathcal{Y} \subset \mathcal{L}$. 
For node $i \in \mathcal{N}$, $v_i$ denotes the magnitude squared voltage.
For node $j \in \mathcal{N^+}$, $p_j$ and $q_j$ denote the net demand of real and reactive power, respectively.
A positive (negative) value of $p_j$ refers to withdrawal (injection); 
similarly for $q_i$.
Net injections at the root node are denoted by $P_0$ and $Q_0$, for real and reactive power, respectively.
These are positive (negative) when power is flowing from (to) the transmission system.
For each line $ij$, with resistance $r_{ij}$ and reactance $x_{ij}$, $l_{ij}$ denotes the magnitude squared current, $P_{ij}$ and $Q_{ij}$ the sending-end real and reactive power flow, respectively.

The \emph{branch} \emph{flow} (AC power flow) equations, are listed below, where we introduce the time index $t$ omitted previously for brevity;
unless otherwise mentioned, $j \in \mathcal{N^+}$, and $t \in \mathcal{T^+}$, with $\mathcal{T} = \{0,1,...,T \}$, $\mathcal{T^+} \equiv \mathcal{T} \setminus \{ 0 \}$, and $T$ the length of the optimization horizon.
\begin{equation} \label{EqSubstation}
P_{0,t} = P_{01,t}  \rightarrow (\lambda_{0,t}^P), \quad Q_{0,t} = Q_{01,t}  \rightarrow(\lambda_{0,t}^Q), \quad \forall t,
\end{equation}
\begin{equation} \label{EqRealBalance}
P_{ij,t} - r_{ij} l_{ij,t} = \sum_{k: j \to k} P_{jk,t} + p_{j,t} \rightarrow (\lambda_{j,t}^P), \quad \forall j,t,
\end{equation}
\begin{equation} \label{EqReactiveBalance}
Q_{ij,t} - x_{ij} l_{ij,t} = \sum_{k: j \to k} Q_{jk,t} + q_{j,t} \rightarrow (\lambda_{j,t}^Q), \quad \forall j,t,
\end{equation}
\begin{equation} \label{EqVoltageDef}
v_{j,t} = v_{i,t} - 2 r_{ij} P_{ij,t} - 2 x_{ij} Q_{ij,t} + \left( r_{ij}^2 + x_{ij}^2 \right) l_{ij,t}, \,\,\, \forall j,t,
\end{equation}
\begin{equation} \label{EqCurrentDef}
 v_{i,t} l_{ij,t} = P_{ij,t}^2 + Q_{ij,t}^2 \quad \forall j,t.
\end{equation}
Briefly, (\ref{EqSubstation})--(\ref{EqReactiveBalance}) define the real and reactive power balance, (\ref{EqVoltageDef}) the voltage drop, and (\ref{EqCurrentDef}) the apparent power but can be also viewed as the definition of the current.
We supplement the model with voltage and current limits, as follows: 
\begin{equation} \label{EqVoltageLimits}
\underline{v}_i \leq v_{i,t} \leq \bar{v}_i \rightarrow ( \underline{\mu}_{i,t}, \bar \mu_{i,t}), \quad \forall i,t,
\end{equation}
\begin{equation} \label{EqCurrentLimits}
l_{ij,t} \leq \bar{l}_{ij}  \rightarrow ( \bar{\nu}_{j,t} ), \quad  \forall j,t,
\end{equation}
where $\underline{v}_i$, $\bar{v}_i$, and $\bar{l}_{ij}$ are the lower voltage, upper voltage, and line ampacity limits (squared), respectively.
Dual variables of constraints (\ref{EqSubstation})--(\ref{EqReactiveBalance}), (\ref{EqVoltageLimits}) and (\ref{EqCurrentLimits}) are shown in parentheses. 

\subsection{DER Models} \label{DERconstraints}

The real (reactive) power net demand $p_{j,t}$ ($q_{j,t}$) includes the aggregate effect of: 
(\emph{i}) conventional demand consumption $p_{d,t}$ ($q_{d,t}$) of load $d \in \mathcal{D}_j$, where $\mathcal{D}_j \subset \mathcal{D} $ is the subset of loads (set $ \mathcal{D}$) connected at node $j$;
(\emph{ii}) consumption $p_{e,t}$ ($q_{e,t}$) of EV $e \in \mathcal{E}_{j,t}$, where $\mathcal{E}_{j,t} \subset \mathcal{E}$ is the subset of EVs (set $\mathcal{E}$) that are connected at node $j$, during time period $t$, and 
(\emph{iii}) generation $p_{s,t}$ ($q_{s,t}$) of PV (rooftop solar) $s \in \mathcal{S}_j$, where $\mathcal{S}_j \subset \mathcal{S}$ is the subset of PVs (set $\mathcal{S}$) connected at node $j$.
For clarity, the definitions of aggregate dependent variables are listed below:
\begin{equation} \label{Pinject}
p_{j,t} = \sum_{d \in \mathcal{D}_j} p_{d,t} + \sum_{e \in \mathcal{E}_{j,t}}  p_{e,t} - \sum_{s \in \mathcal{S}_j} p_{s,t}, \quad\forall j, t,
\end{equation}
\begin{equation} \label{Qinject}
q_{j,t} = \sum_{d \in \mathcal{D}_j} q_{d,t} + \sum_{e \in \mathcal{E}_{j,t}}  q_{e,t} - \sum_{s \in \mathcal{S}_j} q_{s,t} , \quad \forall j, t.
\end{equation}

Next, we provide PV and EV constraints for a general multi-period operational planning problem setting, accommodating smart inverter capabilities and EV mobility (EVs can be connected at different nodes during the time horizon).

\subsubsection{PV Constraints}
Due to the irradiation level $\rho_t$, the PV nameplate capacity $C_s$ is adjusted to $\tilde{C}_{s,t} = \rho_t C_s$, where $\rho_t \in [0,1]$.  
PV constraints ($\forall s \in \mathcal{S}$) are as follows:
\begin{equation} \label{EqPVcon1}
p_{s,t} \leq \tilde{C}_{s,t}, \quad p_{s,t}^2 + q_{s,t}^2 \leq C_s^2, \quad \forall s, t \in \mathcal{T}_I,
\end{equation}
\begin{equation} \label{EqPVcon2}
p_{s,t} = q_{s,t} = 0, \qquad \forall s, t \not \in \mathcal{T}_I,
\end{equation}
with $p_{s,t} \ge 0$, and $\mathcal{T}_I \subset \mathcal{T^+}$ the subset of time periods for which $\rho_t > 0$.
Constraints (\ref{EqPVcon1}) impose limits on real and apparent power (implicitly assuming an appropriately sized inverter), whereas (\ref{EqPVcon2}) imposes zero generation when $\rho_t = 0$.

\subsubsection{EV Constraints}
We consider an EV that is connected for $Z$ intervals, at nodes $j_1,...,j_Z$.
Let $\tau_z^{arr}$ ($\tau_z^{dep}$), for $z = 1,...,Z$, denote the periods at which an arrival (departure) occurs. 
In general, 
intervals $1$ (first) and $Z$ (last) may not entirely fit within the time horizon.
Hence, let $\mathcal{T}^{beg} = \{ \tau_z^{beg} \}$ and $\mathcal{T}^{end} = \{ \tau_z^{end} \}$ be the sets of time periods of interval $z$, for $z=1,...,Z$, denoting an adjusted beginning and end, respectively, considering only the part of the interval within the time horizon.
Also, let $\mathcal{T}_z = \{ \tau_z^{beg}+1,...,\tau_z^{end} \}$ be the set of time periods of interval $z$, for $z = 1,..,Z$, during which the EV is connected at node $j_z$.
Subscript $e$ in the aforementioned sets was omitted for simplicity; it is included next.

The State of Charge (SoC) of EV $e$ is described by variable $u_{e,t}$ for time periods $t \in \mathcal{T}_e^{beg} \cup \mathcal{T}_e^{end}$.
The SoC is reduced by $\Delta u_{z,e}$ after departure $z$ and until arrival $z+1$, for $z = 1,..,Z_e-1$. 
EV constraints ($\forall e \in \mathcal{E}$) are as follows:
\begin{equation} \label{EVCon1}
u_{e,\tau_1^{beg}} = u_e^{init}, \quad \forall e,
\end{equation}
\begin{equation} \label{EVCon2}
u_{e,\tau_z^{end}} = u_{e,\tau_z^{beg}} + \sum_{t \in \mathcal{T}_{e,z}} p_{e,t}, \quad \forall e, z = 1,...,Z_e,
\end{equation}
\begin{equation} \label{EVCon3}
u_{e,\tau_{z+1}^{beg}} = u_{e,\tau_z^{end}} - \Delta u_{e,z}, \quad \forall e, z = 1,...,Z_e-1,
\end{equation}
\begin{equation} \label{EVCon4} 
u_{e,t}^{min} \leq u_{e,t} \leq C_e^B, \quad \forall e, t \in \mathcal{T}_{e}^{end},
\end{equation}
\begin{equation} \label{EVCon5}
p_{e,t}^2 + q_{e,t}^2 \leq C_e^2, \quad 0 \leq p_{e,t} \leq C_{r}, \quad \forall e, t \in { \cup_{z=1}^{Z_e} \mathcal{T}_{e,z} },
\end{equation}
\begin{equation} \label{EVCon6}
p_{e,t} = q_{e,t} = 0, \quad \forall e, t \in {\mathcal{T^+}} \setminus { \cup_{z=1}^{Z_e} \mathcal{T}_{e,z} },
\end{equation}
with $p_{e,t}, u_{e,t} \geq 0$.
Eq. \eqref{EVCon1} initializes the SoC ($u_e^{init}$) at $\tau_1^{beg}$, \eqref{EVCon2} and \eqref{EVCon3} define the SoC at the end/beginning of an interval, after charging/traveling, respectively. 
Constraints \eqref{EVCon4} impose a minimum SoC, $u_{e,t}^{min}$, at the end of an interval as well as the limit of the EV battery capacity, $C_e^B$, whereas \eqref{EVCon5} impose the limits of the charger, $C_e$, (related to the size of the inverter) and the charging rate, $C_{r}$, (related to the capacity of the EV battery charger).
Lastly, \eqref{EVCon6} imposes zero consumption when the EV is not plugged in.

\subsection{Transformer Degradation Model} \label{TransformerModel}

For a given HST of the winding, $\theta^H$, IEEE and IEC Standards \cite{IEEE2011, IEC2018} provide the following exponential representation for the aging acceleration factor, $F_{AA}$:
\begin{equation} \label{FaaExp}
F_{AA} = \exp \left( \frac{15000}{383}-\frac{15000}{\theta^H + 273} \right).
\end{equation}
Both standards provide detailed formulas on the transformer thermal response that have been widely used in simulation studies, but have not been employed in an AC OPF model.
In what follows, we provide the first AC OPF formulation that captures transformer thermal response dynamics taking advantage of the features of the branch flow model.

\begin{figure}[tb]
\centering
\includegraphics[width=3in]{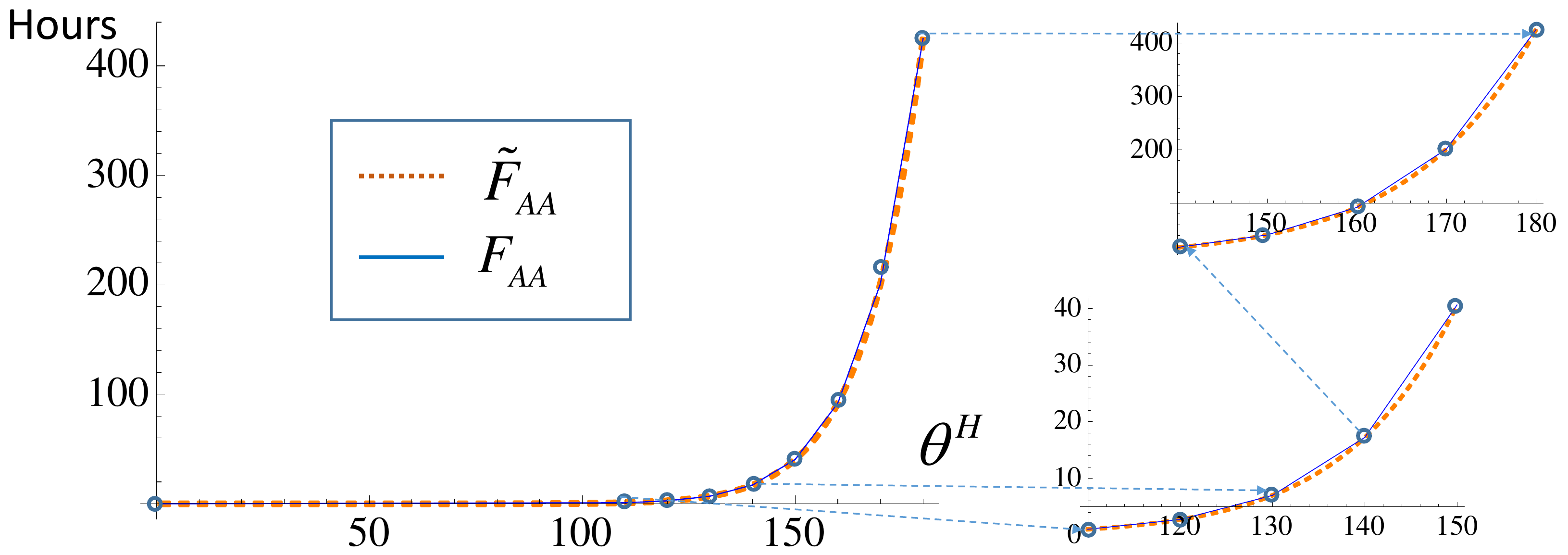}
\caption{Aging acceleration factor $F_{AA}$ vs. HST $\theta^H$ (in $^o$C).
Piecewise linear approximation $\tilde{F}_{AA}$; $M=8$; indicative intervals for $\theta^H$: 0-110, 110-120,...,170-180.}
\label{figFaa}
\end{figure}

First, we consider a piecewise linear approximation of \eqref{FaaExp} denoted by $\tilde{F}_{AA}$ and illustrated in Fig. \ref{figFaa}, which is given by
\begin{equation} \label{Eq:PieceWise} 
\tilde{F}_{AA} = a_{\kappa} \theta^H - b_{\kappa}, \,\, \theta_{\kappa-1}^H \leq \theta^H < \theta_{\kappa}^H, \,\, \kappa = 1,...,M,
\end{equation}
where $M$ is the number of segments.
Introducing indices $y$ and $t$, and substituting $\tilde{F}_{AA}$ with $f_{y,t}$ to simplify the notation, the transformer degradation cost in the optimization horizon is represented by $\sum_{y,t} c_y f_{y,t}$, where $c_y$ is the hourly cost of transformer $y$.
This term (cost) in the objective function, allows us to replace \eqref{Eq:PieceWise} with the following set of inequalities:
\begin{equation} \label{PiecewiseIneq}
f_{y,t} \geq a_{\kappa} \theta^H_{y,t} - b_{\kappa}, \qquad \forall y,t,\kappa,
\end{equation}
with $f_{y,t} \geq 0$. 
Since $c_y > 0$, it is straightforward to show that at least one of the above inequalities should be binding (at equality)\footnote{
Assume that none of the inequalities was binding. 
Then, this solution would not be optimal, as the transformer cost could have been reduced by decreasing the value of $f_{y,t}$ until at least one of the constraints becomes binding.
In fact, at most two constraints can be binding simultaneously (two in the case when the HST is located at a breakpoint of the piecewise linear approximation).}.
Hence, the relaxation of \eqref{Eq:PieceWise} --- which would otherwise require the introduction of binary variables --- to the linear inequality constraints \eqref{PiecewiseIneq} is exact.

Next, we define the HST using linear recursive equations that fit nicely with the branch flow model. 
The HST of transformer $y$ at time period $t$, $\theta^H_{y,t}$, is given by:
\begin{equation} \label{EqThetaH}
\theta^H_{y,t} = \theta^A_{y,t} + \Delta \theta^{TO}_{y,t} + \Delta \theta^H_{y,t} = \theta^{TO}_{y,t} + \Delta \theta^H_{y,t},
\end{equation}
where
$\theta^A_{y,t}$ is the ambient temperature at the transformer location, $\theta^{TO}_{y,t}$ is the top-oil (TO) temperature,
$\Delta \theta^{TO}_{y,t}$ is the TO temperature rise over $\theta^A_{y,t}$,
and $\Delta \theta^H_{y,t}$ is the winding HST rise over $\theta^{TO}_{y,t}$.
The detailed derivations are based on the heat transfer differential equations.
For ease of exposition, they are listed in Appendix \ref{AppA}, summarized in formulas \eqref{Eq:thetaTOapprox} and \eqref{Eq:DthetaHapprox}.

Considering the context of the operational planning problem, we highlight two important characteristics of our model, which greatly facilitate the accurate representation of the transformer thermal response. 
The first one is related to the granularity of the model (indicatively hourly). 
Both $\theta^{TO}_{y,t}$ and $\Delta \theta^H_{y,t}$ are characterized by an oil time constant and a winding time constant, respectively, whose typical values of about 3 hours for the oil and of about 4 minutes for the winding, allow us to employ difference equations at the hourly timescale for the oil whereas assume a steady state for the winding.
Notably, this would still hold if we allowed up to a 15-min granularity of the operational planning problem.
The second characteristic is related to the branch flow model and its decision variable $l_{y,t}$, representing the magnitude squared of the current.
This fact allows us to define the square of the ratio of the transformer load to the rated load, which appears in the HST calculations (see Appendix \ref{AppA}), using variable $l_{y,t}$ and the transformer nominal current (at rated load) squared, denoted by $l_y^N$.
Hence, the equations that define the HST, which we embed in the AC OPF model, are linear. 

Combining \eqref{EqThetaH} and \eqref{Eq:DthetaHapprox}, we define $\theta_{y,t}^{H}$ as a linear equation with variables the TO temperature at time period $t$, $\theta_{y,t}^{TO}$, and the transformer load represented by $l_{y,t}$. 
Substituting $\theta_{y,t}^{TO}$ with $h_{y,t}$, to simplify the notation, and replacing $\theta_{y,t}^H$ in \eqref{PiecewiseIneq}, we get ($\forall y \in \mathcal{Y}, t \in \mathcal{T^+}, \kappa = 1,...,M$): 
\begin{equation} \label{Xf1}
f_{y,t} \geq \alpha_{\kappa} h_{y,t} + \beta_{y, \kappa} l_{y,t} + \gamma_{y,\kappa}
\rightarrow (\xi_{y,t,\kappa}), \, \forall y, t, \kappa,
\end{equation}
with $f_{y,t} \ge 0$, and $\xi_{y,t,\kappa}$ the dual variable.
The coefficients $\alpha_\kappa$, $\beta_{y,\kappa},$ and $\gamma_{y,\kappa}$ can be obtained directly from \eqref{PiecewiseIneq} and \eqref{Eq:DthetaHapprox}. Using the recommended values for the general parameters, they are given by:
\begin{equation} \label{Coef1}
\alpha_\kappa = a_{\kappa}, \,\, \beta_{y,\kappa} =  a_{\kappa} \frac{4 \Delta \bar{\theta}_{y}^H}{ 5 l_y^N}, 
\,\, \gamma_{y,\kappa} = a_{\kappa} \frac{\Delta \bar\theta_{y}^H}{5} - b_{\kappa},
\end{equation}  
where $\Delta \bar \theta_{y}^{H}$ is the rise of HST over TO temperature at rated load, $a_{\kappa}$ and $b_{\kappa}$ are the slope and intercept, respectively, of the $\kappa$-th segment of the piecewise linearization in \eqref{Eq:PieceWise}.

The TO temperature at time period $t$, $\theta_{y,t}^{TO}$, is in turn defined by a linear recursive equation that includes the TO temperature at time period $t - 1$, $\theta_{y,t-1}^{TO}$, and $l_{y,t}$. Indeed, there is also an impact of the ambient temperature, $\theta_{y,t}^{A}$, but this is in fact an input parameter.
Substituting $\theta_{y,t}^{TO}$ with $h_{y,t}$ in \eqref{Eq:thetaTOapprox}, we get the following recursive equation for the TO temperature:
\begin{equation} \label{Xf2}
h_{y,t} = \delta h_{y,t-1} + \epsilon_y l_{y,t} + \zeta_{y,t},\quad \forall y, t.
\end{equation}
The coefficients $\delta$, $\epsilon_{y}$, and $\zeta_{y,t}$ can be obtained by \eqref{Eq:thetaTOapprox}. 
Using again the recommended values, they become:
\begin{equation} \label{Coef2}
\delta = \frac{3}{4}, \,
\epsilon_{y} =  \frac{ R_y \Delta \bar \theta_{y}^{TO} }{5(1 + R_y) l_y^N }, \,
\zeta_{y,t} = \frac{(5 + R_y)\Delta \bar \theta_{y}^{TO}}{20(1 + R_y)} + \frac{\theta_{y,t}^{A}}{4},
\end{equation} 
where $\Delta \bar \theta_{y}^{TO}$ is the rise of TO temperature over ambient temperature at rated load, and $R_y$ is the ratio of load losses at rated load to the losses at no load. Notably, the value of $\delta$ that is always less than $1$ plays an important role in the intertemporal impact of the transformer degradation cost, which will become evident in the next section.\footnote{
From \eqref{Eq:thetaTOapprox}, $\delta = \frac{k_1 \tau^{TO}}{k_1 \tau^{TO}+ \Delta t} < 1$, where $k_1 = 1$, $\tau^{TO} = $3h, and $\Delta t$ the granularity of the optimization problem. 
Hence, for time periods $\Delta t$ equal to 1 hour, 30 minutes or 15 minutes, $\delta$ will be 0.75, 0.857 or 0.923, respectively.
}

Summarizing, the transformer degradation model includes the transformer degradation cost in the objective function, represented by $\sum_{y,t} c_y f_{y,t}$, and linear constraints \eqref{Xf1}, \eqref{Xf2}, with coefficients defined by \eqref{Coef1} and \eqref{Coef2}.

One issue that naturally arises in the context of the detailed and enhanced operational planning problem that internalizes intertemporal transformer cost relations is the impact of decisions at period $t$ that extend beyond the optimization horizon.
Potential remedies and some additional discussion follow.

First, rolling horizon approaches applied in scheduling problems that typically involve future uncertainty are relevant in modeling such intertemporal impacts. 
Second, an extended horizon for the transformer degradation cost might offer a reasonable alternative solution. 
Such an approach should include the next half day, to capture 4 oil time constants, and could assume an estimate for the extended horizon, while accounting for anticipated ambient temperature trend.
A third approach, which is quite suitable for simulation purposes, is to set the following constraint:
\begin{equation} \label{Cycle}
h_{y,T} = h_{y,0} \to (\rho_y), \qquad \forall y, 
\end{equation}
which essentially models the daily 24-hour ahead problem as a cycle repeating over identical days.
This constraint's dual variable, $\rho_y$, captures the future impact on future transformer LoL of loading towards the end of the day. 
Alternatively, we could require some explicit condition for $h_{y,T}$, e.g., equal or less than or equal to some target for the initial condition of the next day, and add a cost in the objective function (pre-calculated by offline studies) for deviating from this target.

\subsection{Optimization Problem Summary} \label{ProblemSummary}

The objective function of the operational planning optimization problem --- referred to as \textbf{Full-opt} --- aims at minimizing the aggregate real and reactive power cost, with $c_t^P$ ($c_t^Q$) denoting the cost for real (reactive) power at the substation, as well as the transformer degradation cost. 
The real power cost $c_t^P$ is typically the LMP at the T\&D interface, whereas $c_t^Q$ can be viewed as the opportunity cost for the provision of reactive power. Hence, Full-opt is defined as follows:\\
\textbf{Full-opt}:
\begin{equation} \label{Obj2}
\text{minimize}
\overbrace{\sum_{t} c^P_t P_{0,t}}^{ \text{Real Power Cost} }
+ \overbrace{\sum_{t} c^Q_t Q_{0,t}}^{\text{ Reactive Power Cost} }
+ \overbrace{\sum_{y,t}{ c_y f_{y,t} }}^{ \text{Transformer Cost} },
\end{equation}
\emph{subject to:} network constraints \eqref{EqSubstation}--\eqref{Qinject}, transformer constraints \eqref{Xf1}, \eqref{Xf2}, \eqref{Cycle}, PV constraints \eqref{EqPVcon1}, \eqref{EqPVcon2}, and  EV constraints \eqref{EVCon1}--\eqref{EVCon6}, with variables $v_{i,t}, l_{ij,t}, f_{y,t}, p_{s,t}, p_{e,t},  u_{e,t}$ nonnegative, and $P_{0,t}, Q_{0,t}, P_{ij,t}, Q_{ij,t}, p_{j,t}, q_{j,t}, h_{y,t}, q_{s,t}, q_{e,t}$ unrestricted in sign.

It is important to note that equality constraint (\ref{EqCurrentDef}) is non-linear and as such introduces a non-convexity to the allowable decision set.
Following \cite{FarivarLow2013}, we relax the equality to an inequality constraint, and substitute (\ref{EqCurrentDef}) with   
\begin{equation} \label{EqCurrentDefRelaxed}
 v_{i,t} l_{ij,t} \geq P_{ij,t}^2 + Q_{ij,t}^2 \quad \forall j,t.
\end{equation}
The resulting relaxed AC OPF problem is a convex SOCP problem, which can be solved efficiently using commercially available solvers.

\section{DLMC Components} \label{DLMCs}

In this section, we provide a rigorous analysis and interpretation of the DLMC components, with emphasis on the new transformer degradation component.
DLMCs represent the dynamic marginal cost for delivering real and reactive power, denoted by P-DLMC and Q-DLMC, respectively, at a specific location and time period.
In the context of the Full-opt operational planning problem, they are obtained by the dual variables of constraints \eqref{EqSubstation}--\eqref{EqReactiveBalance}.
In fact, it is trivial to show that, at the root node, $\lambda_{0,t}^P = c_t^P$, and $\lambda_{0,t}^Q = c_t^Q$.
Hence, in what follows, we will refer to DLMCs as the dual variables of constraints \eqref{EqRealBalance} and \eqref{EqReactiveBalance}.

We note that DLMCs can be defined for any DER schedule reflected in net demand variables $p_{j,t}$ and $q_{j,t}$, $\forall j, t$ (implying, unless otherwise mentioned, $j \in \mathcal{N}^+, t \in \mathcal{T}^+$).
Given these variables, the solution of the AC OPF problem defines an operating point for the distribution network, which, in fact, can be obtained by the solution of the power flow equations \eqref{EqRealBalance}--\eqref{EqCurrentDef}, when the root node voltage is given.
The operating point is described by the real/reactive power flows, voltages and currents, i.e., variables $P_{ij,t}$, $Q_{ij,t}$, $v_{j,t}$, and $l_{ij,t}$.

Let us consider the DLMCs at node $j'$ and time period $t'$, $\lambda_{j',t'}^P$ and $\lambda_{j',t'}^Q$.
The decomposition of the DLMCs employs sensitivity analysis, duality and optimality conditions of the enhanced AC OPF problem.
The sensitivity of the power flow solution \emph{w.r.t.} net demand for real and reactive power, is reflected in the partial derivatives $\frac{\partial P_{ij,t'}}{\partial p_{j',t'}} $, $\frac{\partial Q_{ij,t'}}{\partial p_{j',t'}} $, $\frac{\partial v_{j,t'}}{\partial p_{j',t'}}$, $\frac{\partial l_{ij,t'}}{\partial p_{j',t'}}$, and $\frac{\partial P_{ij,t'}}{\partial q_{j',t'}} $, $\frac{\partial Q_{ij,t'}}{\partial q_{j',t'}} $, $\frac{\partial v_{j,t'}}{\partial q_{j',t'}}$, $\frac{\partial l_{ij,t'}}{\partial q_{j',t'}}$, respectively.
Their calculation requires the solution of a linear system, which is listed in \cite[Appendix A]{PartII}.
More precisely, the calculation of sensitivities \emph{w.r.t.} real net demand at node $j'$ and time period $t'$ (in total $4N$ sensitivities/unknowns) involves the solution of a linear system derived from the power flow equations (\ref{EqRealBalance})-(\ref{EqCurrentDef}), after taking the partial derivatives \emph{w.r.t.} $p_{j',t'}$, at time period $t'$ , for all nodes, thus $4N$ equations.
Similarly for the the sensitivities \emph{w.r.t.} reactive net demand.
Hence, in total, it requires the solution of $2NT$ linear systems (referring to the set of sensitivities \emph{w.r.t.} real and reactive power net demand, calculated at $N$ nodes, and for $T$ time periods), with each system involving $4N$ equations and an equal number of unknowns.
As we noted in Subsection \ref{Contribution}, the ability to rely on the solution of linear systems rather than on brute force perturbation, removes computational and numerical approximation error concerns associated with the need to solve multiple perturbed AC OPF problems.
Furthermore, as noted in \cite{Papavasiliou_2017} the calculations are amenable to parallelization; in our case, each linear system can be solved in parallel.

In what follows, we provide and discuss the unbundled DLMCs to additive components:
\begin{equation} \label{pDLMC}
\begin{split}
\lambda_{j',t'}^P 
= c_{t'}^P 
+ \overbrace{c_{t'}^P  \sum_j r_{ij}\frac{\partial l_{ij,t'}}{\partial p_{j',t'}}}^{\text{Real Power Marginal Losses}}
+ \overbrace{c_{t'}^Q  \sum_j x_{ij}\frac{\partial l_{ij,t'}}{\partial p_{j',t'}}}^{\text{Reactive Power Marginal Losses}} \\
 + \overbrace{\sum_j{ \mu_{j,t'} \frac{\partial v_{j,t'}}{\partial p_{j',t'}}}}^{\text{Voltage Congestion}} 
+ \overbrace{\sum_j { \nu_{j,t'} \frac{\partial l_{ij,t'}}{\partial p_{j',t'}}}}^{\text{Ampacity Congestion}}
+ \overbrace{\sum_{y}{\pi_{y,t'} \frac{\partial l_{y,t'}}{\partial p_{j',t'}}}}^{\text{Transformer Degradation}},
\end{split}
\end{equation}
\begin{equation} \label{qDLMC}
\begin{split}
\lambda_{j',t'}^Q  
= c_{t'}^Q 
+ \overbrace{c_{t'}^P  \sum_j r_{ij}\frac{\partial l_{ij,t'}}{\partial q_{j',t'}}}^{\text{Real Power Marginal Losses}}
+ \overbrace{c_{t'}^Q  \sum_j x_{ij}\frac{\partial l_{ij,t'}}{\partial q_{j',t'}}}^{\text{Reactive Power Marginal Losses}} \\
 + \overbrace{\sum_j{ \mu_{j,t'} \frac{\partial v_{j,t'}}{\partial q_{j',t'}}}}^{\text{Voltage Congestion}} 
+ \overbrace{\sum_j { \nu_{j,t'} \frac{\partial l_{ij,t'}}{\partial q_{j',t'}}}}^{\text{Ampacity Congestion}}
+ \overbrace{\sum_{y}{\pi_{y,t'} \frac{\partial l_{y,t'}}{\partial q_{j',t'}}}}^{\text{Transformer Degradation}},
\end{split}
\end{equation}
where parameter $\pi_{y,t'}$, which is defined and explained next, includes the intertemporal impact of the transformer degradation component.
For both P-DLMC/Q-DLMC, the first component is the marginal cost of real/reactive power at the root node (substation).
The second and the third components represent the contribution of real/reactive power marginal losses. 
Both real/reactive power marginal losses are sensitive to changes in net real/reactive power demand.
The first three terms are obtained by associating equations \eqref{EqRealBalance} and \eqref{EqReactiveBalance} recursively to \eqref{EqSubstation}, and taking the partial derivatives.
This yields:
\begin{equation*} 
\frac{\partial P_{0,t'}}{\partial p_{j',t'}}
=
1 + \sum_{j} r_{ij}\frac{\partial l_{ij,t'}}{\partial p_{j',t'}},
\, \frac{\partial Q_{0,t'}}{\partial p_{j',t'}}
=
\sum_{j} x_{ij}\frac{\partial l_{ij,t'}}{\partial p_{j',t'}} ,
\end{equation*}
\begin{equation*} 
\frac{\partial P_{0,t'}}{\partial q_{j',t'}}
=
\sum_{j} r_{ij}\frac{\partial l_{ij,t'}}{\partial q_{j',t'}},
\,
\frac{\partial Q_{0,t'}}{\partial q_{j',t'}}
=
1 + \sum_{j} x_{ij}\frac{\partial l_{ij,t'}}{\partial q_{j',t'}}.
\end{equation*}
The fourth component reflects voltage congestion, with $\mu_{j,t'} = \bar \mu_{j,t'} - \underline \mu_{j,t'}$.
It is nonzero only at nodes with a binding voltage constraint, i.e., when $\bar \mu_{j,t'} > 0$ or $\underline \mu_{j,t'} > 0$ (obviously they cannot be positive simultaneously).
The fifth component represents ampacity congestion.
It is non zero only at lines with a binding ampacity constraint, i.e., when $\nu_{j,t'} > 0$.
The voltage and ampacity congestion components can be directly derived by appending the active (binding) constraints in the Full-opt Lagrangian. 
Alternatively, they can be obtained by the optimality conditions. 
Last, and perhaps most important, the sixth component represents the impact on transformer degradation costs that are intertemporally coupled with real and reactive power injections during preceding time periods. 
We next derive the formula for $\pi_{y,t'}$ and elaborate on it.

The derivation is identical for P-DLMC and Q-DLMC (note that parameter $\pi_{y,t'}$ is the same for both), hence, without loss of generality, we employ the P-DLMC component.
Applying equation \eqref{Xf2} recursively for $h_{y,t}$:
\begin{equation*}
h_{y,t} = \delta^t h_{y,0} + \sum_{\tau = 1}^t \delta^{t-\tau} \left( \epsilon_y l_{y,\tau} + \zeta_{y,\tau} \right), \quad \forall y, t,
\end{equation*}
and taking the partial derivatives of $h_{y,t}$ \emph{w.r.t.} $p_{j',t'}$, we get:
\begin{equation} \label{hytRecursive2}
\frac{\partial h_{y,t}}{\partial p_{j',t'}} = \epsilon_y \delta^{t-t'}  \frac{\partial l_{y,t'}}{\partial p_{j',t'}}  ,\quad \forall y, t \geq t'.
\end{equation}
From the binding transformer constraints \eqref{Xf1}, using optimality conditions and \eqref{hytRecursive2}, we obtain:
\begin{equation} \label{Margfyt}
\begin{split}
\sum_{y,t}{c_y \frac{\partial f_{y,t}}{\partial p_{j',t'}}} 
= \sum_{y,t, \kappa}{\xi_{y,t,\kappa} \left( \alpha_{\kappa} \frac{\partial h_{y,t}}{\partial p_{j',t'}} + \beta_{y, \kappa} \frac{\partial l_{y,t}}{\partial p_{j',t'}}  \right)} \qquad \quad\\
= \left[ \epsilon_y \sum_{t=t'}^T { \delta^{t-t'} \left( \sum_{\kappa}{\xi_{y,t,\kappa} \alpha_{\kappa}} \right) } 
+ \sum_{\kappa} \xi_{y,t',\kappa} \beta_{y,\kappa} \right] \frac{\partial l_{y,t}}{\partial p_{j',t'}},   
\end{split}
\end{equation}
where the sum of dual variables $\sum_{\kappa}\xi_{y,t,\kappa}$ should equal $c_y$.\footnote{
Obviously, the sum $\sum_{\kappa} \xi_{y,t,\kappa}$ involves only the binding constraints (otherwise $\xi_{y,t,\kappa} = 0$), and may contain up to two terms, in case the solution is found at a breakpoint. 
If only one constraint is active, then the respective dual $\xi_{y,t,\kappa}$ should equal $c_y$.  
}
Considering also constraint \eqref{Cycle} and appending it in the Lagrangian, we obtain the following term:
\begin{equation} \label{CycleMarginal}
\rho_y \frac{\partial h_{y,T}}{\partial p_{j',t'}} = \epsilon_y \delta^{T-t'}  \rho_y \frac{\partial l_{y,t'}}{\partial p_{j',t'}},
\end{equation}
which denotes the impact that extends beyond the optimization horizon.
Combining \eqref{Margfyt}, which denotes the impact within the optimization horizon, and \eqref{CycleMarginal}, $\pi_{y,t'}$ is given by:
\begin{equation} \label{piyt}
\pi_{y,t'} =  \epsilon_y \left( \sum_{t=t'}^T  \delta^{t-t'} \tilde \alpha_{y,t} +  \delta^{T-t'}  \rho_y \right)
+ \eta_y \tilde \alpha_{y,t'},
\end{equation}
where, using \eqref{Coef1}, $\tilde \alpha_{y,t} = \sum_{\kappa} \xi_{y,t,\kappa} a_{\kappa}$, $\eta_y = \frac{4 \Delta \bar{\theta}_{y}^H}{ 5 l_y^N}$, and $\epsilon_y$ is given by \eqref{Coef2}.

Notably, $\pi_{y,t'}$ involves a summation over time that is related to the transformer thermal response dynamics and captures subsequent time period costs.
The first term in \eqref{piyt} that is multiplied by $\epsilon_y$ refers to the contribution of the load at time period $t'$ to the TO temperature at $t'$ and subsequent time periods as well as the impact that extends beyond the optimization horizon (involving dual variable $\rho_y$), whereas the second term that is multiplied by $\eta_y$ refers to the contribution to the winding temperature rise over the TO.
Recall that the winding time constant is small compared to our time period length of one hour, hence its impact is limited to time period $t'$, whereas the oil time constant is larger, and hence the intertemporal impact on subsequent time periods is not negligible.
This intertemporal impact is discounted by a factor $\delta^{t-t'}$ (where $\delta = 3/4 < 1$), but it also considers the slope\footnote{
Recall that this slope is the piecewise linear approximation of the exponential which is large when the HST is high reflecting current and past period loading of the transformer.} 
of transformer LoL, $a_{\kappa}$, and applies a higher weight to overloaded time periods.
The impact that extends beyond the optimization horizon is discounted by $\delta^{T-t'}$, implying that this impact may become significant for hours that are closer to the end of the horizon.

\section{Conclusions and Further Research} \label{Conclusions}
We have provided an enhanced AC OPF for granular marginal costing on distribution networks (DLMCs) including asset degradation costs exemplified by transformer marginal LoL.
We have also derived and analyzed additive DLMC components, and illustrated the intertemporal characteristics of the transformer degradation component.
The DLMCs estimated by the proposed AC OPF are compatible with optimal DER self-scheduling that is fully adapted to the optimal distribution network configuration. 
As such, the DLMCs represent the spatiotemporal marginal costs on the distribution network that co-optimize network and DER scheduling decisions. 
Whereas we have considered the real and reactive power marginal cost at the substation bus, or equivalently at the T\&D interface, to be given by setting it equal to the corresponding wholesale market LMP, it is possible that after all distribution networks and their DER schedules have been co-optimized, the Wholesale Transmission market may adjust producing new LMPs.
Although the real power LMPs are most likely to be insensitive to real power DLMCs that affect primarily distribution node injections but less so the aggregate T\&D power exchange, it is possible that reactive power compensation at the T\&D node may be impacted by the distribution network reactive power schedule.
This is a topic of future work that may require ISOs to provide reasonable, possibly extended LMP type, price functions at the T\&D interface.

Lastly, our ongoing and future research includes extending our DLMC framework to a 3-phase network representation, and employing decomposition approaches and distributed algorithms \cite{MolzahnEtAl_2017} to deal with multiple feeder systems and large numbers of diverse DERs transcending EVs and PVs and including micro generators, smart buildings with pre-cooling/heating capable HVAC, smart appliances, storage, and new technologies.

\appendices

\section{Transformer HST Calculations} \label{AppA}

Transformer HST is described in \cite{IEC2018} as a function of time, for varying load and ambient temperature, using exponential equations, and difference equations.
Both methods represent solution variation to the heat transfer differential equations.

For the TO temperature, $\theta^{TO}_t$, where we dropped transformer index $y$, the differential equation is given by:
\begin{equation} \label{EqDif}
\Delta \bar \theta^{TO} \left(\frac{1 + K_t^2 R}{1 + R}\right)^n = k_{1} \tau^{TO}\frac{d\theta^{TO}_t}{dt} + \theta^{TO}_t - \theta^A_t,
\end{equation}
where $\Delta \bar \theta^{TO}$ is the rise of TO temperature over ambient temperature at rated load, $K_t$ is the ratio of the (current) load to the rated load, $R$ is the ratio of load losses at rated load to no-load losses, $\tau^{TO}$ is the oil time constant with recommended value 3 hours, $k_{1}$ and $n$ are constants with recommended values 1 and 0.8, respectively.
Since the granularity for a day-ahead problem is less than half the recommended value of $\tau^{TO}$ ($\Delta t = 1$ hour), we can employ the difference equations and get the following recursive formula:
\begin{equation} \label{Eq:thetaTO}
\begin{split}
\theta_t^{TO} & = \frac{k_{1} \tau^{TO}}{k_{1} \tau^{TO} + \Delta t} \theta_{t-1}^{TO} \\
 & + \frac{\Delta t}{k_{1} \tau^{TO} + \Delta t} \left[ \Delta \bar \theta^{TO} \left(  \frac{1+K_t^2 R}{1+R}\right)^n  + \theta_t^{A} \right].
\end{split}
\end{equation}

For the winding HST rise over $\theta^{TO}_t$, $\Delta \theta^H_t$, it can be shown that both \cite{IEEE2011} and \cite{IEC2018} yield the same results for distribution (small) transformers.
Because the winding time constant $\tau^w$ has an indicative value of about 4 min (much less than the hourly granularity), the transient behavior $1-\exp \left(-\frac{\Delta t}{\tau^w} \right)  \approx 1$, vanishes.
Hence, using \cite{IEEE2011}, we get:
\begin{equation} \label{Eq:DthetaH}
\Delta \theta_t^{H} = \Delta \bar \theta^{H} \left( K_t^2\right) ^{m},
\end{equation}
where $\Delta \bar \theta^{H}$ is the rise of HST over TO temperature at rated load, and $m$ is a constant with recommended value 0.8.

The load ratio $K_{t}$ can be defined \emph{w.r.t.} the transformer nominal current (at rated load), denoted by $I^N$.
Using variables $l_{t}$ (omitting the transformer index $y$), we have $K_{t}^2 = {l_{t}}/{l^N}$, where $l^N = ({I^N})^2$.
Hence, approximating the terms $\left(\frac{1+K_t^2 R}{1+R}\right)^n$ and $\left( K_t^2\right) ^{m}$ in (\ref{Eq:thetaTO}) and (\ref{Eq:DthetaH}), respectively, using the 1st order Taylor expansion of $K_{t}^2$ around 1, (equivalently of $l_t$ around $l^N$), and replacing $K_{t}^2$ by ${l_{t}}/{l^N}$, we get:
\begin{equation} \label{EqApprox1}
\left(\frac{1+{K_t}^2 R}{1+R} \right)^n \approx  
 \frac{nR}{(1+R)} \frac{l_t}{l^N} + \frac{1+(1-n)R}{1+R},
\end{equation}
\begin{equation} \label{EqApprox2}
\left( K_t^2\right) ^{m}  \approx m  \frac{l_t}{l^N} + 1 - m .
\end{equation}

Using \eqref{EqApprox1} and \eqref{EqApprox2}, \eqref{Eq:thetaTO} and \eqref{Eq:DthetaH} yield:
\begin{equation} \label{Eq:thetaTOapprox}
\begin{split}
\theta_t^{TO} = & \frac{k_1 \tau^{TO}}{k_1 \tau^{TO}+ \Delta t} \theta_{t-1}^{TO}\\ 
& + \frac{\Delta t}{k_{1} \tau^{TO} + \Delta t}
    \frac{n R}{1 + R} \Delta \bar \theta^{TO} \frac{l_t}{l^N} \\
& + \frac{\Delta t}{k_{1} \tau^{TO} + \Delta t} 
    \left[
    \frac{ 1 + (1-n)R }{1+R} \Delta \bar \theta^{TO}
    + \theta_t^{A} \right],
\end{split} 
\end{equation}
\begin{equation} \label{Eq:DthetaHapprox}
\Delta \theta_t^{H} =  m \Delta \bar \theta^{H} \frac{l_t}{l^N}  + (1-m) \Delta \bar \theta^{H}.
\end{equation}

Lastly, we note that the initial value $\theta_0^{TO}$ --- in case it is not known --- can be obtained assuming a steady state, i.e., setting the derivative in the differential equation (\ref{EqDif}) to zero, and replacing $K_0^2$ by $l_0/l^N$, which yields:
\begin{equation} \label{Eq:thetaTO_0}
\theta_0^{TO} = \theta_0^{A} +  \Delta \bar \theta^{TO} \left[\frac{1+ (l_0/l^N) R }{1+R}\right]^n. 
\end{equation}

\ifCLASSOPTIONcaptionsoff 
  \newpage
\fi

\end{document}